%% file: tyler_convexification.tex
\documentclass{article}

\usepackage{spconf,amsfonts,amsmath,cite,graphicx,bm,amssymb,amsthm,enumerate,epsfig,psfrag,cases,mathtools}
\include{ilyas_commands}

\usepackage{fontenc}
\usepackage{inputenc}

\begin{document}

\ninept

\title{Covariance Estimation in Elliptical Models with Convex Structure}

\name{Ilya Soloveychik and Ami Wiesel\thanks{This work was partially supported by the Intel Collaboration Research Institute for Computational Intelligence and Kaete Klausner Scholarship.
}}
\address{The Rachel and Selim Benin School of Computer Science and Engineering,\\ The Hebrew University of Jeursalem, Israel}

\maketitle

\begin{abstract}
We address structured covariance estimation in Elliptical distribution. We assume it is a priori known that the covariance belongs to a given convex set, e.g., the set of Toeplitz or banded matrices. We consider the General Method of Moments (GMM) optimization subject to these convex constraints. Unfortunately, GMM is still non-convex due to objective. Instead, we propose COCA - a convex relaxation which can be efficiently solved. We prove that the relaxation is tight in the unconstrained case for a finite number of samples, and in the constrained case asymptotically. We then illustrate the advantages of COCA in synthetic simulations with structured Compound Gaussian distributions. In these examples, COCA outperforms competing methods as Tyler's estimate and its projection onto a convex set. 
\end{abstract}


\section{Introduction}
Covariance matrix estimation is a fundamental problem in the field of statistical signal
processing. Many other algorithms for detection and inference rely on accurate covariance estimates \cite{krim1996two,dougherty2005research}. The problem is well understood in the Gaussian unstructured case. But becomes significantly harder when the underlying distribution is non-Gaussian, for example in Elliptical distributions, and when there is prior knowledge on the structure. In this paper, we propose a unified framework for covariance estimation in Elliptical distributions with general convex structure.

Over the last years there was a great interest in covariance estimation with known structure. The motivation to these works is that in many modern applications the dimension of the underlying distribution is large and there are not enough samples to estimate it correctly. The prior information on the structure reduces the degrees of freedom in the model and allows accurate estimation with a small number of samples. This is clearly true when the structure is exact, but also when it is approximate due to the well known bias-variance tradeoff. Prior knowledge on the structure can originate from the physics of the underlying phenomena, e.g., \cite{fuhrmann1991application, roberts2000hidden, pollock2002circulant, stoica2011spice}, or from similar datasets, e.g., adjacent cells in radar systems \cite{wang1994adaptive}. When the structure is defined using a convex set, a natural and computationally efficient solution is to project the naive unstructured estimators onto this set.

Many covariance structures are easily represented in a convex form. Probably the most classical structure is the Toeplitz model. It arises naturally in the analysis of stationary time series which are used in a wide range of applications including radar imaging, target detection, speech recognition, and communication systems, \cite{snyder1989use, fuhrmann1991application, roberts2000hidden}. Toeplitz matrices are also used
to model the correlation of cyclostationary processes in periodic time series \cite{dahlhaus1989efficient}. In other settings the number of parameters can be reduced by assuming that the covariance matrix sparse\cite{cai2012optimal,bickel2008regularized}.  A popular sparse model is the banded covariance, which is associated with time-varying moving average models\cite{bickel2008regularized}. Another important example of a convex structure is the SPICE estimator, which was proposed in \cite{stoica2011spice} to treat high-dimensional arrays processing problems, where the covariance structure is approximated by a low-dimensional linear combination of known rank one matrices. Finally, a structure that recently attracted considerable attention involves low-rank matrices. In the last decade, all of these structures have been successfully considered in the Gaussian case.

In a different line of works, there is an increasing interest in robust covariance estimation for non-Gaussian distributions \cite{pascal2008covariance}. Significant attention is being paid to the family of Elliptical distributions, which includes as particular cases generalized Gaussian distribution, compound Gaussian and many other \cite{frahm2004generalized}.  Elliptical models are commonly used to measure radar clutter \cite{conte1991modelling}, noise and interference in indoor and outdoor mobile communication channels \cite{middleton1973man} and other applications. For these purposes, robust covariance estimators were developed \cite{maronna1976robust}. In particular, \cite{tyler1987distribution} proposed a robust scatter estimator which has become widely used \cite{abramovich2007time, pascal2008covariance,bandiera2010knowledge}. One of the most prominent disadvantage of these methods is that they involve non-convex optimization problem, thus making imposition of additional constraints rather difficult. One of the options to cure this obstacle is geodesic convexity. It has been recently shown that some of the popular M-estimators are in fact g-convex, which significantly simplifies their treatment \cite{wiesel2012geodesic,soloveychik2013group}. But still, if one wants to impose an additional constraint on the scatter matrix it must be formulated in a form of a g-convex set, rather than a classical convex set. 

In the present work we derive COCA -  COnvexly ConstrAined Covariance Matching estimator. It is based on the principle of Generalized Method of Moments (GMM) \cite{matyas1999generalized}. It searches for a covariance of a given convex structure that minimizes the norm of a simple moment's identity. This identity is in fact the optimality condition of Tyler's estimate. COCA tries to simultaneously satisfy this condition while constraining the structure. Unfortunately, this requires the solution to a high dimensional non-convex minimization. Instead, we propose a relaxation and express COCA as a standard convex optimization
with linear matrix inequalities which can be computed using off-the-shelf numerical solvers. Interestingly, we prove two promising results. First, in the unconstrained case, COCA is tight and identical to Tyler's estimate. This result basically ``convexifies'' Tyler's estimate. Second, in the structured case, COCA is asymptotically tight and hence consistent. Finally, we demonstrate the finite sample advantages of COCA over existing methods using synthetic numerical simulations.

The paper is organized in the following way. First, we formulate the problem and briefly describe the existing solutions: the sample covariance, Tyler's estimator and the projection method. We then introduce GMM, derive its convex relaxation denoted by COCA and show that its coincides with Tyler's estimator in the unconstrained case. After this we prove that adding convex structure does not affect asymptotic consistency. Finally, we provide numerical examples and applications demonstrating the performance advantages of COCA.

We denote by $\mathcal{P}(p)$ the closed cone of symmetric positive semi-definite $p \times p$ matrices. 

\section{Problem Formulation}
Consider a $p$ dimensional, zero mean, Elliptically distributed random vector $\x$. Such a vector can be represented as \cite{frahm2004generalized}
\begin{equation}
\x = r {\bf{\Lambda}} \u,
\end{equation}
where $\u$ is a $k$ dimensional random vector, uniformly distributed on the unit hypersphere, $r$ is a nonnegative random variable, $\Lambda \in \mathbb{R}^{p \times k}$ \cite{frahm2004generalized}.  The random variable $r$ is called the generating variate of $\x$ and it is stochastically independent of $\u$. We assume that the distribution of this random variable is unknown.

The parameter $\C^{\text{True}} = {\bf{\Lambda}} {\bf{\Lambda}}^T$ is referred to as the dispersion or shape matrix of $\x$ and coincides with its covariance matrix (up to a scaling factor). In many applications, it is common to assume prior information on the structure of this matrix. In particular, we assume that it belongs to a known convex subset $\mathcal{S} \subset \mathcal{P}(p)$. Typical examples of such subsets are:

\begin{itemize}

\item {\bf{Toeplitz}}: In stationary time series, the covariance between the $i$-th and the $j$-th components depend only on the the difference $|i-j|$. Such kind of processes is encountered very often in many engineering areas including statistical signal processing, radar imaging, target detection, speech recognition, and communications systems, \cite{burg1982estimation,snyder1989use, fuhrmann1991application, roberts2000hidden, wiesel2013time,dahlhaus1989efficient}.

\item {\bf{Banded}}: A natural approach to covariance modeling is to formulate the reduction in statistical relation using the notion of independence or correlation, which corresponds to sparsity in the covariance matrix \cite{bickel2008regularized}. Assuming that $i$-th element of the random vector is uncorrelated with the $j$-th if $|i-j|>k$ leads to $k$-banded structure, also known as time varying moving average models.

\item {\bf{Low rank}}: One of the most common covariance models involves a low dimensional principal subspace plus white noise \cite{ganz1990convergence}. A typical convex representation of such models is $\C^{\text{True}}=\X+\sigma^2\I$ together with a bound on the nuclear norm of $\X\in \mathcal{P}(p)$. In this model, $\sigma^2$ is the known and fixed variance of the noise.

\item {\bf{Linear parameterization}}: Many interesting models can be expressed as a linear combination of known matrices. In particular, a modern approach to estimation of direction of arrivals of multiple signals involves a covariance of the form $\C^{\text{True}} = \sum_{i=1}^k  p_i \a_i\a_i^T$ where $\a_i$ constitute a dense grid of possible directions, and $p_i$ are their corresponding coefficients. Typically, the $l_1$ norm of these sparse coefficients is constrained. See for example \cite{stoica2011spice}.

\end{itemize}

We can now state the problem addressed in this paper. Let $\x_i, i=1,\dots,n$ be independent and identically distributed (i.i.d) copies of $\x$ with $\C^{\text{True}} \in \mathcal{S}$. Given these realizations and knowledge of $\mathcal{S}$, we are interested in estimation of the matrix $\C^{\text{True}}$.

\section{Existing solutions}
\subsection{Sample Covariance}
The classical solution to the above covariance estimation problem is the sample covariance matrix defined by
\begin{equation}
\C^{\text{Sample}} = \frac{1}{n} \sum_{i=1}^n \x_i \x_i^T.
\label{c_sam}
\end{equation}
The sample covariance estimator always exists and is asymptotically consistent in any distribution with bounded moments by the Law of Large Numbers. In the Gaussian case when $n\geq p$, it also maximizes the likelihood and is asymptotically efficient. In the non-Gaussian case, it has been extensively studied \cite{vershynin2012close}, any is generally suboptimal. Furthermore, it does not exploit any additional structure knowledge.
\subsection{Tyler's M-estimation approach}
The most popular approach to covariance estimation in Elliptical distribution is due to Tyler \cite{tyler1987distribution}. This estimator is defined as the fixed point solution to:
\begin{equation}
\C^{\text{Tyler}} = \frac{p}{n} \sum_{i=1}^n \frac{\x_i \x_i^T}{\x_i^T {\left[\C^{\text{Tyler}}\right]}^{-1} \x_i}.
\label{tyler_formula}
\end{equation}
Since $\C^{\text{Tyler}}$ is defined only up to scale, it has to be fixed by some additional constraint like $\Tr {\C^{\text{Tyler}}}=1$. When $n>p$, it has been proven that a simple fixed point iteration converges to this unique solution \cite{maronna1976robust}. This estimator is asymptotically consistent in all Elliptical distributions. In fact, it is has been shown to maximize the likelihood of the normalized samples
\begin{equation}
\s=\frac{\x}{||\x||_2} = \frac{r \Lambda \u}{||r \Lambda \u||_2} = \frac{\Lambda \u}{||\Lambda \u||_2},
\label{unit_ge}
\end{equation}
which are independent of the values of the generating variates. The advantages of Tyler's estimator are its simplicity and robustness. Its drawbacks are that it does not exist if $n<p$ and does not exploit known structure. In \cite{bandiera2010knowledge} knowledge based variants of the fixed point iteration were proposed without convergence analysis. Recently, regularized and structured versions of Tyler's estimate were proposed in \cite{abramovich2007time,chen2011robust,wiesel2011regularized,wiesel2012geodesic,wiesel2012geodesic,soloveychik2013group} based on the theories of concave Perron Frobenius and geodesic convexity. Unfortunately, these approaches are limited in its their modeling capabilities are cannot deal with general convex models as described above.

\subsection{Projection}
A reasonable approach to explore the covariance structure is to use a projection. Given any estimator $\hat\C$, e.g., the sample covariance or Tyler, its projection onto ${\mathcal{S}}$ is defined as
\begin{equation}
\C^{\text{Proj}}_{\mathcal{S}}\(\hat\C\)= \argmin_{\M \in \mathcal{S}}||\M-\hat\C||, \label{proj_est}
\end{equation}
where $\|\cdot\|$ is some norm. For simple structures as described above, the projection is a convex optimization problem which can be efficiently solved using standard numerical packages, e.g., CVX, \cite{cvx, gb08}. The main advantage is that, when $\C^{\text{True}} \in \mathcal{S}$, the projection $\C^{\text{Proj}}$ is usually closer to $\C^{\text{True}}$ than $\hat\C$ is. The disadvantage is that it requires a two-step solution which does not take into couple the distribution properties and the structure information and is therefore suboptimal.

\section{COCA-estimator}
In this section, we propose COCA - the COnvexly ConstrAined covariance estimator for Elliptical distributions. Unlike the existing solutions, COCA exploits both the Elliptical nature and the structure of the underlying distribution. COCA is based on the Generalized Method of Moments \cite{matyas1999generalized} together with an asymptotically tight convex relaxation.

 The underlying principle behind COCA is the following identity  \cite{frahm2004generalized, frahm2010generalization}:
\begin{equation}
\E \left(p \frac{\x_i \x_i^T}{\x_i^T \[\C^{\text{True}}\]^{-1} \x_i}\right) = \C^{\text{True}},
\label{main_lem}
\end{equation}
Indeed, Tyler's estimator is just the sample based solution that satisfies this identity. When there is an insufficient number of samples and a constraint of the structure, such a solution does not necessarily exist. Instead, we propose the Generalized Method of Moments \cite{matyas1999generalized} which seeks an approximate solution to
\begin{equation}
\begin{aligned}
& \underset{\C\in\mathcal{S}_1}{\text{min}}
& & \left|\left|\C-\frac{p}{n}\sum_{i=1}^n \frac{\x_i \x_i^T}{\x_i^T \C^{-1} \x_i} \right|\right|, \\
\label{gmm_pr}
\end{aligned}
\end{equation}
where $||\cdot||$ is some norm and we ensure uniqueness by defining
\begin{equation}
\mathcal{S}_1 = \{\M \in \mathcal{S}|\Tr \M =1\}. \label{main_tr}
\end{equation}
Intuitively, this optimization tries to simultaneously solve Tyler and project it on the prior structure. By choosing an adaptive weighted norm, an optimal solution to (\ref{gmm_pr}) would result in an asymptotically consistent and efficient estimator \cite{matyas1999generalized,ottersten1998covariance}). Unfortunately, the objective is non-convex and it is not clear how to find its global solution in a tractable manner.

In what follows, we propose a convex relaxation of (\ref{gmm_pr}) that allows a computationally efficient solution. First, let us introduce the auxiliary variables $d_i,i=1,\dots,n$:
\begin{equation}
\begin{aligned}
& \underset{\C\in\mathcal{S}_1,d_i}{\text{min}}
& & \left|\left|\C-\frac{1}{n}\sum_{i=1}^n d_i \x_i \x_i^T \right|\right| \\
& \text{subject to}
& & d_i = \frac{p}{\x_i^T \C^{-1} \x_i}, i=1 \dots n.
\end{aligned}
\end{equation}
This problem is not convex due to the equality constraints. We suggest to relax them to inequalities:
\begin{equation}
\begin{aligned}
& \underset{\C\in\mathcal{S}_1,d_i}{\text{min}}
& & \left|\left|\C-\frac{1}{n}\sum_{i=1}^n d_i \x_i \x_i^T \right|\right| \\
& \text{subject to}
& & d_i \leq \frac{p}{\x_i^T \C^{-1} \x_i}, i=1 \dots n, \\
& & & d_i \geq 0, i=1 \dots n.
\end{aligned}
\end{equation}
This relaxed problem is a convex minimization. To see this, it is instructive to use Schur's complement formulas and express the inequalities $d_i \leq \frac{p}{\x_i^T \C^{-1} \x_i}, i=1 \dots n$ as convex linear matrix inequalities (LMI):
\begin{equation}
\C^{\text{COCA}} = \arg \left\{
\begin{aligned}
& \underset{\C\in\mathcal{S}_1,d_i}{\text{min}}
& & \left|\left|\C-\frac{1}{n}\sum_{i=1}^n d_i \x_i \x_i^T \right|\right| \\
& \text{subject to}
& & \C \succeq \frac{1}{p} d_i \x_i \x_i^T, \forall i=1 \dots n,\\
& & & d_i \geq 0, \forall i=1 \dots n.
\end{aligned}
\right.
\label{sdp_f}
\end{equation}
COCA can be efficiently computed by standard semi-definite programming solvers, e.g., CVX, \cite{cvx, gb08}.

The non-relaxed version of COCA in (\ref{gmm_pr}) can be considered optimal in many ways. The interesting question is how tight is the relaxation. We now provide two promising results in this direction.
\begin{theorem} \label{theorem1}
In the unstructured case $\mathcal{S} = \mathcal{P}(p)$ with $n \geq p$, COCA is unique up to a positive scaling factor and coincides with Tyler's estimator.
\begin{proof}
It is known that when $n \geq p$, (\ref{sdp_f}) has at least one solution which results in a zero objective value. It is Tyler's estimator which satisfies
\begin{equation}
d_i = \frac{p}{\x_i^T \C^{-1} \x_i}, i = 1 \dots n. \nonumber
\end{equation}

It remains to show that there are no other feasible solutions which result in a zero objective. Indeed, assume in contradiction that there is such a solution, and that $\C=\frac{1}{n} \sum_{i=1}^n d_i \x_i \x_i^T$. Multiply each inequality $d_i \leq \frac{p}{\x_i^T \C^{-1} \x_i}$ by the matrix $\x_i \x_i^T$ for $i=1 \dots n$ and sum up to obtain
\begin{equation}
\C = \frac{1}{n} \sum_{i=1}^n d_i \x_i \x_i^T \preceq \frac{p}{n} \sum_{i=1}^n \frac{\x_i \x_i^T}{\x_i^T \C^{-1} \x_i} = f(\C). \label{tyl_ineq}
\end{equation}

The inequality (\ref{tyl_ineq}) reads now as $\C \preceq f(\C)$.
As stated in the Corollary V.I from \cite{pascal2008covariance}, this implies that $\C$ is the fixed point of $f$: $\C = f(\C)$ (See also \cite{maronna1976robust}), which is exactly the definition of Tyler's estimator in (\ref{tyler_formula}). Thus proving that it is the only solution to (\ref{sdp_f}) up to a positive scaling factor.
\end{proof}
\end{theorem}

\begin{theorem} \label{theorem2}
In the structured case, COCA is an asymptotically consistent estimator of the true shape matrix $\C^{\text{True}} \in \mathcal{S}$.
\begin{proof}
Due to space limitations, we defer this technical proof to the journal version of this paper. In brief, the idea is that, as $n \rightarrow \infty$, choosing $\C=\C^{\text{Tyler}}\in \mathcal{S}$ and $d_i=\frac{p}{\x_i^T\[\C^{\text{Tyler}}\]^{-1}\x_i}$, yields an asymptotic zero objective. Then we show that the inequality (\ref{tyl_ineq}) also holds asymptotically and therefore this solution is unique (up to a scaling factor).
\end{proof}
\end{theorem}

\section{Numerical results}
In this section, we demonstrate the advantages of COCA using numerical simulations. For simplicity, we chose the norm in (\ref{sdp_f}) and the norm of the projector operator (\ref{proj_est}) as the spectral norm. We investigated the performance benefits of COCA when the true shape matrix was either Toeplitz or banded. We compared the following estimators: $\C^{\text{Sample}}$ in (\ref{c_sam}), $\C^{\text{Tyler}}$ in (\ref{tyler_formula}), $\C^{\text{Proj}}$ in (\ref{proj_est}) and $\C^{\text{COCA}}$ in (\ref{sdp_f}). In $\C^{\text{Proj}}$ we projected Tyler's estimator when it existed and the sample covariance otherwise.

For each number of samples $n$ we generated $1000$ sets of independent, Elliptically distributed $20$-dimensional samples and calculated the empirical MSE for all the estimators. The samples were generated as compound Gaussian  $\x = \sqrt{\tau} \mathbf{v}$, where the random variable $\tau \sim \chi^2$ and the random vector $\mathbf{v}$ was zero-mean normally distributed with covariance matrix $\C^{\text{True}}$.

\subsection{Toeplitz Covariance Matrix}
The $20 \times 20$ Toeplitz shape matrix was obtained as $\C^{\text{True}} = \F \D \F^T$, where $\F$ is the $20$-dimensional DFT matrix and $\D$ is a diagonal matrix with eigenvalues $1,\dots,20$, e.g., \cite{pollock2002circulant}. This matrix is complex Hermitian, and all the theory developed above applies by replacing the transpose operators with conjugate transposes. The average results are reported in Fig. 1. It is easy to see the performance advantage of COCA over all previously known estimators over a wide range of number of samples.

\begin{figure}
 \center\includegraphics[scale=0.6]{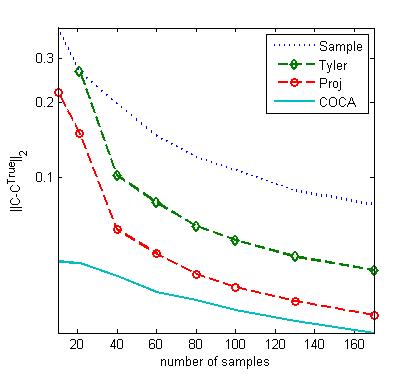}\caption{COCA in Toeplitz models.}
\end{figure}

\subsection{Banded Covariance Matrix}
For the true banded covariance matrix we took a symmetric matrix with the numbers $21,\dots,40$ on the diagonal, $1,\dots,19$ and $1,\dots,18$ on the first and second sub-diagonals respectively. The averaged errors are reported in Fig. 2. Here too, COCA outperforms the previous approaches.

\begin{figure}
 \center\includegraphics[scale=0.6]{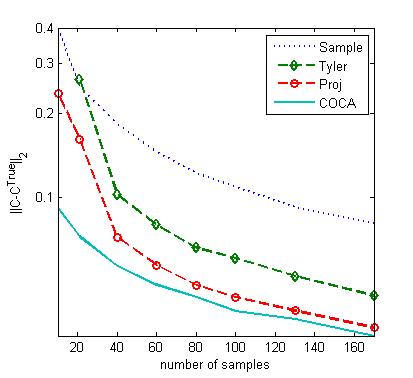}\caption{COCA in Banded models}
\end{figure}

\section{Discussion}
We proposed a novel COCA estimator for structured covariance estimation in Elliptical distributions. In two examples, we demonstrate that it is more accurate than Tyler's classical estimator and its projection. The most important benefits of COCA are that it is defined via a convex program, thus admitting any convex structure; and that it exists when $n < p$, in which case Tyler's estimators does not exist.

There are several ways of extending this result. First, it can be extended to a general M-estimator and not only the particular case of Tyler's estimator. Second, weighted and adaptive norms should be considered for ensuring statistical efficiency as explained in \cite{matyas1999generalized,ottersten1998covariance}. Finally, our current analysis only guarantees asymptotic consistency and future work will address its non-asymptotic properties in comparison to existing methods.

\end{document}

%% file: ilyas_commands.tex
\renewcommand{\(}{\left(}
\renewcommand{\)}{\right)}
\renewcommand{\[}{\left[}
\renewcommand{\]}{\right]}

\renewcommand{\a}{\mathbf{a}}

\newcommand{\s}{\mathbf{s}}

\newcommand{\E}{\mathbf{E}}

\newcommand{\D}{\mathbf{D}}

\newcommand{\F}{\mathbf{F}}
\newcommand{\x}{\mathbf{x}}

\newcommand{\I}{\mathbf{I}}

\newcommand{\C}{\mathbf{C}}

\newcommand{\M}{\mathbf{M}}

\renewcommand{\u}{\mathbf{u}}

\newcommand{\X}{\mathbf{X}}

\newcommand{\Tr}[1]{{\rm{Tr}}\left(#1\right)}

\newcommand{\End}[1]{{\rm{End}}}

\renewcommand{\arg}[1]{{\rm{arg}}#1}

\newtheorem{theorem}{Theorem}

\DeclareMathOperator*{\argmin}{\arg\!\min}

%% file: tyler_convexification.bbl
\begin{thebibliography}{10}

\bibitem{krim1996two}
H.~Krim and M.~Viberg.
\newblock Two decades of array signal processing research: the parametric
  approach.
\newblock {\em IEEE Signal Processing Magazine}, 13(4):67--94, 1996.

\bibitem{dougherty2005research}
E.~R. Dougherty, A.~Datta, and C.~Sima.
\newblock Research issues in genomic signal processing.
\newblock {\em IEEE Signal Processing Magazine}, 22(6):46--68, 2005.

\bibitem{fuhrmann1991application}
D.~R. Fuhrmann.
\newblock Application of {T}oeplitz covariance estimation to adaptive
  beamforming and detection.
\newblock {\em IEEE Transactions on Signal Processing}, 39(10):2194--2198,
  1991.

\bibitem{roberts2000hidden}
W.~J. Roberts and Y.~Ephraim.
\newblock Hidden {M}arkov modeling of speech using {T}oeplitz covariance
  matrices.
\newblock {\em Speech Communication}, 31(1):1--14, 2000.

\bibitem{pollock2002circulant}
D.~S. Pollock.
\newblock Circulant matrices and time-series analysis.
\newblock {\em International Journal of Mathematical Education in Science and
  Technology}, 33(2):213--230, 2002.

\bibitem{stoica2011spice}
P.~Stoica, P.~Babu, and J.~Li.
\newblock {SPICE}: A sparse covariance-based estimation method for array
  processing.
\newblock {\em IEEE Transactions on Signal Processing}, 59(2):629--638, 2011.

\bibitem{wang1994adaptive}
H.~Wang and L.~Cai.
\newblock On adaptive spatial-temporal processing for airborne surveillance
  radar systems.
\newblock {\em IEEE Transactions on Aerospace and Electronic Systems},
  30(3):660--670, 1994.

\bibitem{snyder1989use}
D.~L. Snyder, J.~A. O'Sullivan, and M.~I. Miller.
\newblock The use of maximum likelihood estimation for forming images of
  diffuse radar targets from delay-{D}oppler data.
\newblock {\em IEEE Transactions on Information Theory}, 35(3):536--548, 1989.

\bibitem{dahlhaus1989efficient}
R.~Dahlhaus.
\newblock Efficient parameter estimation for self-similar processes.
\newblock {\em The Annals of Statistics}, pages 1749--1766, 1989.

\bibitem{cai2012optimal}
T.~T. Cai, Zh. Ren, and H.~H. Zhou.
\newblock Optimal rates of convergence for estimating {T}oeplitz covariance
  matrices.
\newblock {\em Probability Theory and Related Fields}, pages 1--43, 2012.

\bibitem{bickel2008regularized}
P.~J. Bickel and E.~Levina.
\newblock Regularized estimation of large covariance matrices.
\newblock {\em The Annals of Statistics}, pages 199--227, 2008.

\bibitem{pascal2008covariance}
F.~Pascal, Y.~Chitour, J.~P. Ovarlez, P.~Forster, and P.~Larzabal.
\newblock Covariance structure maximum-likelihood estimates in compound
  {G}aussian noise: Existence and algorithm analysis.
\newblock {\em IEEE Transactions on Signal Processing}, 56(1):34--48, 2008.

\bibitem{frahm2004generalized}
G.~Frahm.
\newblock Generalized elliptical distributions: theory and applications.
\newblock {\em Universit{\"a}t zu K{\"o}ln}, 2004.

\bibitem{conte1991modelling}
E.~Conte, M.~Longo, and M.~Lops.
\newblock Modelling and simulation of non-{R}ayleigh radar clutter.
\newblock {\em IEE Proceedings F on Radar and Signal Processing},
  138(2):121--130, 1991.

\bibitem{middleton1973man}
D.~Middleton.
\newblock Man-made noise in urban environments and transportation systems:
  Models and measurements.
\newblock {\em IEEE Transactions on Communications}, 21(11):1232--1241, 1973.

\bibitem{maronna1976robust}
R.~A. Maronna.
\newblock Robust {M}-estimators of multivariate location and scatter.
\newblock {\em The annals of statistics}, pages 51--67, 1976.

\bibitem{tyler1987distribution}
D.~E. Tyler.
\newblock A distribution-free {M}-estimator of multivariate scatter.
\newblock {\em The Annals of Statistics}, 15(1):234--251, 1987.

\bibitem{abramovich2007time}
Y.~I. Abramovich, N.~K. Spencer, and M.~D. Turley.
\newblock Time-varying autoregressive ({TVAR}) models for multiple radar
  observations.
\newblock {\em IEEE Transactions on Signal Processing}, 55(4):1298--1311, 2007.

\bibitem{bandiera2010knowledge}
F.~Bandiera, O.~Besson, and G.~Ricci.
\newblock Knowledge-aided covariance matrix estimation and adaptive detection
  in compound-{G}aussian noise.
\newblock {\em IEEE Transactions on Signal Processing}, 58(10):5391--5396,
  2010.

\bibitem{wiesel2012geodesic}
A.~Wiesel.
\newblock Geodesic convexity and covariance estimation.
\newblock {\em IEEE Transactions on Signal Processing}, 60(12):6182--6189,
  2012.

\bibitem{soloveychik2013group}
I.~Soloveychik and A.~Wiesel.
\newblock Group symmetry and non-{G}aussian covariance estimation.
\newblock {\em arXiv preprint arXiv:1306.4103}, 2013.

\bibitem{matyas1999generalized}
L.~M{\'a}ty{\'a}s.
\newblock Generalized method of moments estimation.
\newblock {\em Cambridge University Press}, 5, 1999.

\bibitem{burg1982estimation}
J.~P. Burg, D.~G. Luenberger, and D.~L. Wenger.
\newblock Estimation of structured covariance matrices.
\newblock {\em Proceedings of the IEEE}, 70(9):963--974, 1982.

\bibitem{wiesel2013time}
A.~Wiesel, O.~Bibi, and A.~Globerson.
\newblock Time varying autoregressive moving average models for covariance
  estimation.
\newblock {\em IEEE Transactions on Signal Processing}, 61(11):2791--2801,
  2013.

\bibitem{ganz1990convergence}
M.~W. Ganz, R.~Moses, and S.~Wilson.
\newblock Convergence of the {SMI} and the diagonally loaded {SMI} algorithms
  with weak interference [adaptive array].
\newblock {\em IEEE Transactions on Antennas and Propagation}, 38(3):394--399,
  1990.

\bibitem{vershynin2012close}
R.~Vershynin.
\newblock How close is the sample covariance matrix to the actual covariance
  matrix?
\newblock {\em Journal of Theoretical Probability}, 25(3):655--686, 2012.

\bibitem{chen2011robust}
Y.~Chen, A.~Wiesel, and A.~O. Hero.
\newblock Robust shrinkage estimation of high-dimensional covariance matrices.
\newblock {\em IEEE Transactions on Signal Processing}, 59(9):4097--4107, 2011.

\bibitem{wiesel2011regularized}
A.~Wiesel.
\newblock Regularized covariance estimation in scaled {G}aussian models.
\newblock {\em 4th IEEE International Workshop on Computational Advances in
  Multi-Sensor Adaptive Processing (CAMSAP)}, pages 309--312, 2011.

\bibitem{cvx}
M.~Grant and S.~Boyd.
\newblock {CVX}: Matlab software for disciplined convex programming, version
  2.0 beta.
\newblock http://cvxr.com/cvx, September 2013.

\bibitem{gb08}
M.~Grant and S.~Boyd.
\newblock Graph implementations for nonsmooth convex programs.
\newblock In V.~Blondel, S.~Boyd, and H.~Kimura, editors, {\em Recent Advances
  in Learning and Control}, Lecture Notes in Control and Information Sciences,
  pages 95--110. Springer-Verlag Limited, 2008.

\bibitem{frahm2010generalization}
G.~Frahm and U.~Jaekel.
\newblock A generalization of {T}yler’s {M}-estimators to the case of
  incomplete data.
\newblock {\em Computational Statistics \& Data Analysis}, 54(2):374--393,
  2010.

\bibitem{ottersten1998covariance}
B.~Ottersten, P.~Stoica, and R.~Roy.
\newblock Covariance matching estimation techniques for array signal processing
  applications.
\newblock {\em Digital Signal Processing}, 8(3):185--210, 1998.

\end{thebibliography}
